\documentclass[reqno,11pt]{amsart}
\usepackage[bottom=3cm, right=3cm, left=3cm, top=3cm]{geometry}
\title[Sign Patterns in a Two Colored Partition Companion series]{Sign Patterns in a Two Colored Partition Companion series}
\usepackage{amssymb,amsmath,amsthm,epsfig,graphics,latexsym}
\usepackage{mathrsfs}
\usepackage[hidelinks]{hyperref}
\theoremstyle{definition}

 \usepackage{tikz}
\usetikzlibrary{arrows.meta,positioning,shapes.geometric}
\theoremstyle{plain}
\newtheorem{theorem}    {Theorem}

\newtheorem{lemma}      {Lemma}
\newtheorem{corollary}  {Corollary}
\newtheorem{conjecture} {Conjecture}
\theoremstyle{remark}
\newtheorem{example}{{\bf Example}}
\newtheorem{remark}{{\bf Remark}}
\numberwithin{equation}{section}
 \newcommand{\Mod}[1]{\ (\mathrm{mod}\ #1)}

\newmuskip\pFqskip
\mathchardef\pFcomma=\mathcode`, % keep a copy of the comma

 \newmuskip\pGqskip
\mathchardef\pGcomma=\mathcode`, % keep a copy of the comma

\author[A. Dhar]{Aritram Dhar}
\address{Department of Mathematics, University of Florida, Gainesville, FL 32611, USA}
\email{aritramdhar@ufl.edu, aritramdhar@gmail.com}

\author[A. Goswami]{Ankush Goswami}
\address{School of Mathematical and Statistical Sciences, University of Texas Rio Grande Valley, Edinburg, TX 78541, USA}
\email{ankushgoswami3@gmail.com, ankush.goswami@utrgv.edu}

\author[M. Tripathi]{Mohit Tripathi}
\address{Department of Mathematics and Statistics, Texas Tech University, Lubbock, TX 79410, USA}
\email{mohit.tripathi@ttu.edu}

\date{}

\begin{document}
\begin{abstract}
We study three closely related questions arising from the recent work of Andrews and El Bachraoui on the two-color partition series
\[
S_1(q)=\sum_{n\ge0}s_1(n)q^n=\sum_{a\ge0}q^a(-q^{a+1};q)_\infty^2
\]
and its odd companion, denoted by $T_o(q)$. First, for the eta-normalized companion
\[
C(q)=(q;q)_\infty T_o(q)=\sum_{n\ge0}c(n)q^n,
\]
we prove a strong form of the Andrews--El Bachraoui sign conjecture that $\limsup c(n)=+\infty$ and $\liminf c(n)=-\infty$. Second, we construct an involution using the Franklin-type involution of Chen and Liu to combinatorially explain Andrews--El Bachraoui congruence for $s_1(n)$ modulo 4. Third, we introduce a two-variable refinement of $S_1(q)$ which records the difference of the lengths of the two colored components. Using Ramanujan's partial theta identity and heat-type equations, we obtain a closed formula for all of its even moments, and we show that the eta-normalized moments satisfy termwise congruences modulo $p^{2k}$ for every prime $p\equiv5\pmod6$.
\end{abstract}
\maketitle

% \begin{center}
% \small
% $^{a}$ Department of Mathematics,
% University of Florida,\
% Gainesville, FL 32611, USA\\
% \texttt{aritramdhar@ufl.edu, aritramdhar@gmail.com}

% \vspace{0.5em}

% $^{b}$ School of Mathematical and Statistical Sciences,
% University of Texas Rio Grande Valley,\
% Edinburg, TX 78541, USA\\
% \texttt{ankushgoswami3@gmail.com, ankush.goswami@utrgv.edu}
% % \texttt{runqiao.li@utrgv.edu}

% \vspace{0.5em}

% $^{c}$ Department of Mathematics and Statistics,
% Texas Tech University,\
% Lubbock, TX 79410, USA\\
% \texttt{mohit.tripathi@ttu.edu}

% %\vspace{0.5em}

% %$^{*}$ Corresponding author:
% %\texttt{ankushgoswami3@gmail.com} (\texttt{Ankush Goswami})
% \end{center}

\section{Introduction}\label{s1}
Throughout the paper, we assume $|q|<1$ and use the standard notation
\[
(a;q)_n:=\prod_{j=0}^{n-1}(1-aq^j),
\qquad
(a;q)_\infty:=\prod_{j=0}^{\infty}(1-aq^j),
\]
for the $q$-Pochhammer symbol, together with the shorthand $(a,b;q)_\infty:=(a;q)_\infty(b;q)_\infty$.  We refer to Andrews' text \cite{AndrewsTOP} and to Gasper--Rahman \cite{GasperRahman} for background on $q$-series notation and transformations. Partial theta functions and especially identities of Ramanujan--Andrews type, have long been connected with signed partition identities, false theta phenomena, and Franklin-type cancellations; see, for instance, Alladi, Warnaar, and Andrews--Berndt \cite{Alladi,AndrewsBerndt,Warnaar}.

Recently, Andrews and El Bachraoui \cite{AB} introduced the two-color partition series
\begin{equation}\label{eq:S1intro}
S_1(q):=\sum_{n\ge0}s_1(n)q^n
=\sum_{a\ge0}q^a(-q^{a+1};q)_\infty^2.
\end{equation}
Combinatorially, $s_1(n)$ counts two-color partitions of $n$ into distinct parts in which the smallest part occurs in one prescribed color, while every larger part may occur in either color or in both colors \cite{AB}.  Andrews and El Bachraoui proved the following complete modulo $4$ description.

\begin{theorem}[Andrews--El Bachraoui]\label{thm1AB}
For every $n\ge0$,
\[
s_1(n)\equiv
\begin{cases}
(-1)^r\pmod4, & \text{if } n=r(r+1)/2\text{ for some }r\ge0,\\
0\pmod4, & \text{otherwise.}
\end{cases}
\]
In particular, $s_1(n)$ is odd if and only if $n$ is a triangular number.
\end{theorem}

They also introduced the odd companion
\begin{equation}\label{eq:T_o(q)}
T_o(q):=\sum_{n\ge0}t_o(n)q^n
=\sum_{m\ge0}\frac{q^{2m}(-q;q)_{2m}}{(q;q^2)_{m+1}},
\end{equation}
and the eta-normalized companion
\begin{equation}\label{eq:C(q)}
C(q):=\sum_{n\ge0}c(n)q^n=(q;q)_\infty T_o(q).
\end{equation}
One of the striking features of these series is that $S_1(q)$, $T_o(q)$, and $C(q)$ are tied together by partial-theta decompositions and theta-product normalizations.  In particular, Andrews and El Bachraoui proved the exact relation
\begin{equation}\label{eq:ABcompanion}
S_1(q)=\Theta(q)+4qT_o(q),
\qquad
\Theta(q):=\sum_{r\ge0}(-1)^r q^{r(r+1)/2}.
\end{equation}
Thus, away from triangular indices, the congruence of $s_1(n)$ modulo $8$ is exactly the parity question for $T_o(q)$. Andrews and El Bachraoui also proved the following sign relations for the coefficients of $C(q)$ \cite{AB}.
\begin{theorem}[Andrews--El Bachraoui]\label{thm2AB}
For all $n\ge0$,
\[
c(25n+8)=c(25n+13)=c(25n+18)=c(25n+23)=0,
\]
and
\[
c(25n+28)=-c(n).
\]
In particular, $c(n)$ changes sign infinitely often and assumes both positive and negative values infinitely often.
\end{theorem}

At the end of their paper, Andrews and El Bachraoui proposed the following conjecture \cite{AB}.
\begin{conjecture}[Andrews--El Bachraoui]\label{conj1AB}
We have
\[
\limsup_{n\to\infty}c(n)=+\infty,
\qquad
\liminf_{n\to\infty}c(n)=-\infty.
\]
\end{conjecture}

Our first main result is the following stronger statement.
\begin{theorem}\label{thmstrongconj1AB}
For every integer $k\ge0$, there are infinitely many integers $n$ such that $c(n)=2^k$, and infinitely many integers $n$ such that $c(n)=-2^k$.  In particular,
\[
\limsup_{n\to\infty}c(n)=+\infty,
\qquad
\liminf_{n\to\infty}c(n)=-\infty.
\]
\end{theorem}

The second theme of the paper concerns combinatorial proof of Theorem \ref{thm1AB}, thereby answering Problem 2 in \cite{AB}.
\begin{theorem}
\label{thm:defectfiltration}
Theorem \ref{thm1AB} admits a signed involutive proof.
% : If \(8N+1\) is not represented by the quadratic form \(x^2+2y^2\), then
% \(
% s_1(N)\equiv0\pmod8.
% \)
% Equivalently, if some prime \(p\equiv5,7\pmod8\) divides \(8N+1\) to an odd power, then
% \(
% s_1(N)\equiv0\pmod8.
% \)
\end{theorem}
\begin{remark}
We note here that Andrews and El Bachraoui also conjectured a congruence for $s_1(n)$ modulo 8 \cite[Conjecture 2]{AB}. This has now been proved by Liu--Xia \cite{LiuXia} and, independently, by Banerjee--Bringmann \cite{BanerjeeBringmann}.
\end{remark}

The third theme of the paper is a two-variable refinement of $S_1(q)$.  Define
\begin{equation}
\label{eq:Vdef}
 \mathcal{V}(z;q):=\sum_{a\geq0}q^a
 (-zq^{a+1};q)_\infty(-z^{-1}q^{a+1};q)_\infty
 =\sum_{n\geq0}\sum_{m\in\mathbb{Z}}N(m,n)z^m q^n.
\end{equation}
Combinatorially, \(N(m,n)\) counts triples \((a;\lambda,\mu)\), where \(a\geq0\) and the partitions \(\lambda\) and \(\mu\) have distinct parts larger than \(a\), such that
\[
 a+|\lambda|+|\mu|=n,
 \qquad
 \ell(\lambda)-\ell(\mu)=m.
\]
Here $\ell(\cdot)$ denotes the number of parts.  Clearly $\mathcal V(1;q)=S_1(q)$, so that $s_1(n)=\sum_{m}N(m,n)$, and $\ell(\lambda)-\ell(\mu)$ is a rank-type statistic on the objects counted by $s_1(n)$. Interchanging \(\lambda\) and \(\mu\) gives \(N(m,n)=N(-m,n)\).  Hence all odd moments vanish.  For \(k\geq0\), put
\begin{equation}
\label{eq:momentdef}
 M_{2k}(n):=\sum_{m\in\mathbb{Z}}m^{2k}N(m,n),
 \qquad
 \mathcal{M}_{2k}(q):=\sum_{n\geq0}M_{2k}(n)q^n,
\end{equation}
in analogy with the rank moments of Atkin and Garvan \cite{AtkinGarvan}.  In Section~\ref{s:moments} we show that Ramanujan's partial theta identity splits $\mathcal V(z;q)$ into a theta product times a partial theta function, minus a second partial theta function, and that each piece satisfies a heat-type equation.  This yields a closed formula for every $\mathcal M_{2k}(q)$ (Theorem~\ref{thm:moment}).  For instance, with $\delta:=q\,\frac{d}{dq}$, $q=e^{2\pi i\tau}$, and $E_2(\tau)=1-24\sum_{n\ge1}\sigma_1(n)q^n$,
\[
 \mathcal{M}_2
 =6\delta S_1+\mathcal{E} S_1-2\delta\Theta+(\mathcal{E}-1)\Theta,
 \qquad
 \mathcal{E}=\frac{5E_2(\tau)-8E_2(2\tau)+9}{12}.
\]
Our final main result shows that, after multiplication by $(q;q)_\infty$, the moments are governed by the binary quadratic form $x^2+3y^2$, which is the same form that controls $c(n)$ in Section~\ref{s2}.  Define integers \(d_{2k}(n)\) by
\begin{equation}
\label{eq:ddef}
 (q;q)_\infty\mathcal{M}_{2k}(q)=:\sum_{n\geq0}d_{2k}(n)q^n,
\end{equation}
and let $v_p$ denote the $p$-adic valuation.

\begin{theorem}
\label{thm:cong}
Let \(k\geq1\), and let \(p\equiv5\pmod6\) be prime.  Then
\begin{equation*}
d_{2k}(n)\begin{cases}=0& \text{if }v_p(6n+1)\text{ is odd},\\
\equiv 0\pmod{p^{2ka}}&\text{if }v_p(6n+1)=2a.
\end{cases}
\end{equation*}
In particular,
\begin{equation}
\label{eq:progression}
 d_{2k}\!\left(pn+\frac{5p-1}{6}\right)
 \equiv0\pmod{p^{2k}}.
\end{equation}
\end{theorem}

For example, $d_2(5n+4)\equiv0\pmod{25}$, $d_4(5n+4)\equiv0\pmod{625}$, and $d_2(11n+9)\equiv0\pmod{121}$.  The vanishing part of the argument also applies when $k=0$, and in that case it recovers, and extends to every prime $p\equiv5\pmod6$, the vanishing statements of Theorem~\ref{thm2AB}; see Corollary~\ref{cor:cvanish}.

% The contribution of the present note is not merely to restate this obstruction, but to organize it through a simple $2$-adic filtration.  The identity
% \begin{equation}\label{eq:local_intro}
% (1+q^r)^2=(1-q^r)^2+4q^r
% \end{equation}
% separates $S_1(q)$ into defect layers.  Modulo $4$, only the zero-defect layer remains, and this layer is governed by a Franklin-type involution of Chen and Liu \cite{ChenLiu}.  Modulo $8$, the one-defect layer enters, and the fixed or blocked objects in this layer are controlled by the binary quadratic form $x^2+2y^2$.  Modulo $16$, the same filtration exposes a new secondary obstruction: one needs the one-defect layer modulo $4$, which produces a ternary theta-type correction rather than a purely binary condition.

The paper is organized as follows.  In Section~\ref{s2}, we prove Theorem~\ref{thmstrongconj1AB}.  In Section~\ref{s:combnew}, we give the combinatorial interpretations of the modulo $4$ congruence for $s_1(n)$, thereby proving Theorem \ref{thm:defectfiltration}.  In Section~\ref{s:moments}, we study the refinement $\mathcal V(z;q)$: we establish its partial theta decomposition, prove the moment formula (Theorem~\ref{thm:moment}) and Theorem~\ref{thm:cong}, and record two special values.  We conclude with some open questions in Section~\ref{s:conclusion}.

\section{Proof of Theorem \ref{thmstrongconj1AB}}\label{s2}
We use the signed representation formula for \(c(n)\).  Put $D=6n+7$.
Then \(c(n)\) is a signed count of positive solutions of
\[
x^2+3y^2=4D.
\]
More precisely,
\begin{equation}\label{RD}
c(n)
=
R_{1,3}(D)-R_{5,3}(D)+R_{7,1}(D)-R_{11,1}(D),
\end{equation}
where \(R_{a,b}(D)\) counts positive integer solutions of $x^2+3y^2=4D$
with $x\equiv a\pmod{12},
\,
y\equiv b\pmod4$. Thus the residue classes $(1,3),\,(7,1)$ contribute with plus sign, while $(5,3),\,(11,1)$ contribute with minus sign. We begin with the elementary representation
\[
28=1^2+3\cdot 3^2.
\]
This corresponds to \(D=7\), and the pair $(x,y)=(1,3)$ lies in the plus class
\[
(x,y)\equiv(1,3)\pmod{(12,4)}.
\]
Now choose primes $p_j=u_j^2+3v_j^2$ with $u_j\equiv1\pmod{12},
\,
v_j\equiv0\pmod4$, and with \(v_j/u_j\) very small.  More precisely, choose positive
numbers \(\varepsilon_j\) such that
\[
\sum_{j=1}^{\infty}\varepsilon_j<\eta,
\]
where \(\eta>0\) is chosen small enough so that rotating either of the
two first-quadrant representations of $28=x^2+3y^2$ by any angle of size less than \(\eta\) keeps it in the first quadrant, and rotating any of the other sign choices by an angle of size less than
\(\eta\) keeps it outside the first quadrant.  By Hecke's prime theorem for Grössencharacters, equivalently the prime-ideal theorem in sectors with ray-class conditions
\cite[Ch.~XV, Thm.~6 and Example~2]{LangANT}
(see also Hecke's original formula \cite[p.~38, Eq.~(52)]{Hecke1920}), applied in \(\mathbb Z[\sqrt{-3}]\), every fixed sector
\(
0<\arg(u+\sqrt{-3}v)<\varepsilon
\)
contains infinitely many prime elements, we may
choose infinitely many distinct primes \(p_j\) satisfying the congruences
above and
\(
0<\arctan\left(\frac{\sqrt3\,v_j}{u_j}\right)<\varepsilon_j.
\)
For a fixed \(k\), define
\[
D_k=7p_1p_2\cdots p_k,
\qquad
n_k=\frac{D_k-7}{6}.
\]
Since each \(p_j\equiv1\pmod6\), we have $D_k\equiv7\pmod6$, so \(n_k\) is an integer. The reason these primes are useful is the multiplication identity
\[
(x^2+3y^2)(u^2+3v^2)
=
(xu-3yv)^2+3(xv+yu)^2=(xu+3yv)^2+3(yu-xv)^2,
\]
% and also the conjugate identity
% \[
% (x^2+3y^2)(u^2+3v^2)
% =
% (xu+3yv)^2+3(yu-xv)^2.
% \]
Thus, each time we multiply by one prime \(p_j=u_j^2+3v_j^2\), we get
two possible new representations.  Starting from $28=1^2+3\cdot3^2$
and multiplying successively by $p_1,p_2,\ldots,p_k$, we obtain \(2^k\) representations of
\[
x^2+3y^2=4D_k.
\]
We now check that all these \(2^k\) representations are counted with
plus sign.  Suppose a representation lies in the class
\[
(x,y)\equiv(1,3)\pmod{(12,4)}.
\]
Using
\[
u_j\equiv1\pmod{12},
\qquad
v_j\equiv0\pmod4,
\]
we see from the two multiplication formulas that the new \(x\)-coordinate
is still congruent to \(1\pmod{12}\), and the new \(y\)-coordinate is
still congruent to \(3\pmod4\).  Hence both choices preserve the class $(1,3)$. Therefore all \(2^k\) representations coming from the starting pair $(1,3)$ remain in the plus class $(1,3)$. The small-angle condition ensures that these representations remain
positive.  Indeed, the factor $u_j+\sqrt{-3}\,v_j$
has angle
\[
\arctan\left(\frac{\sqrt3\,v_j}{u_j}\right),
\]
which is smaller than \(\varepsilon_j\).  After \(k\) steps, the total
change in direction is at most
\[
\varepsilon_1+\cdots+\varepsilon_k<\eta.
\]
By the choice of \(\eta\), the resulting points still lie in the first
quadrant.  Thus the corresponding \(x\) and \(y\) remain positive. It remains to explain why the negative classes do not contribute.  The
other first-quadrant representation of
$28=x^2+3y^2$ is $28=5^2+3\cdot1^2$. This pair lies in the class
\[
(x,y)\equiv(5,1)\pmod{(12,4)},
\]
which is not one of the four signed classes appearing in the formula \eqref{RD} for
\(c(n)\).  The same congruence calculation shows that multiplying by any
of the primes \(p_j\) preserves this class.  Hence all representations
generated from \((5,1)\) remain uncounted. Finally, the elementary uniqueness theory for representations by \(x^2+3y^2\) \cite{Cox} says that, since $D_k=7p_1p_2\cdots p_k$ is squarefree and each \(p_j\) is represented by \(u_j^2+3v_j^2\), every
representation of $x^2+3y^2=4D_k$ is obtained from one of the representations of $28=x^2+3y^2$ by making the two choices above at each prime factor, together with
possible sign changes.  By our choice of \(\eta\), the sign-changed
starting points remain outside the first quadrant after the small
rotations, so they do not contribute to the positive-solution count. Consequently, among the four signed classes in the formula for \(c(n)\),
only the plus class $(1,3)$ contributes.  It contributes exactly \(2^k\) times.  Thus
\[
R_{1,3}(D_k)=2^k,
\]
while
\[
R_{5,3}(D_k)=R_{7,1}(D_k)=R_{11,1}(D_k)=0.
\]
Because the primes \(p_j\) can be chosen in infinitely many ways, this
gives infinitely many \(n_k\) with $c(n_k)=2^k$. Hence
\[
\limsup_{n\to\infty}c(n)=+\infty.
\]
The negative direction follows from the known
relation $c(25n+28)=-c(n)$. Applying this to the sequence \(n_k\), we get
\[
c(25n_k+28)=-c(n_k)=-2^k.
\]
Since \(k\) is arbitrary, this gives
\[
\liminf_{n\to\infty}c(n)=-\infty.
\]
The result follows.\qed

\section{Proof of Theorem \ref{thm:defectfiltration}}\label{s:combnew}
In this section we give a signed interpretation of the modulo $4$ congruence for $s_1(n)$.  The key point is the following simple identity:
\begin{equation}\label{eq:localmod4new}
(1+q^r)^2=(1-q^r)^2+4q^r.
\end{equation}
Modulo $4$, $S_1(q)$ and its signed generating function $F(q)$ (defined below) are congruent. By constructing an involution on the signed object using a Franklin-type involution of Chen and Liu \cite{ChenLiu}, we show that the fixed points of our involution are triangular-number staircases, proving Theorem \ref{thm1AB}. Using \eqref{eq:S1intro} and \eqref{eq:localmod4new}, we get
\begin{equation}\label{eq:S1tosignednew}
S_1(q)\equiv F(q):=\sum_{a\ge0}q^a(q^{a+1};q)_\infty^2\pmod4.
\end{equation}
The right-hand side is the signed generating function for triples $X=(a;\lambda,\mu)$, where $a\ge0$ and $\lambda,\mu$ are distinct-part partitions with all parts strictly larger than $a$.  We assign
\[
w(X)=a+|\lambda|+|\mu|,
\qquad
\operatorname{sgn}(X)=(-1)^{\ell(\lambda)+\ell(\mu)}.
\]
Thus
\begin{equation}\label{eq:signedFnew}
F(q)=\sum_X \operatorname{sgn}(X)q^{w(X)}.
\end{equation}
These are exactly the triples counted by $N(m,n)$ in \eqref{eq:Vdef}, and since $\operatorname{sgn}(X)=(-1)^{\ell(\lambda)-\ell(\mu)}$, we have $F(q)=\mathcal V(-1;q)$.

We now describe a complete sign-reversing involution on these triples.  Rather than forcing a long case-by-case map directly on $(a;\lambda,\mu)$, we transfer the problem to the distinct-part model used by Chen and Liu in their Franklin-type involution for squares \cite{ChenLiu}.  Define
\begin{equation}\label{eq:Phi}
\Phi(a;\lambda,\mu)=\rho:=\{2a\}\cup\{2r:r\in\lambda\}\cup\{2s-1:s\in\mu\}.
\end{equation}
Then $\rho$ is a partition into distinct nonnegative parts whose smallest part is even.  If $o(\rho)$ denotes the number of odd parts of $\rho$, then
\begin{equation}\label{eq:weightPhi}
|\rho|+o(\rho)=2w(a;\lambda,\mu),
\end{equation}
and
\begin{equation}\label{eq:signPhi}
(-1)^{\ell(\rho)-1}=(-1)^{\ell(\lambda)+\ell(\mu)}.
\end{equation}
Therefore $F(q)$ is the weighted signed generating function of distinct nonnegative partitions with smallest part even.

Chen and Liu \cite[Theorem 5.1]{ChenLiu} constructed an explicit Franklin-type involution on precisely such distinct-part partitions $\rho$. Denote their involution by $\psi$.  It preserves the statistics $|\rho|$, and their proof confirms that the involution also preserves $o(\rho)$, and hence, $w(a;\lambda,\mu)=(|\rho|+o(\rho))/2$; reverses the parity of $\ell(\rho)-1$, and has fixed points
\begin{equation}\label{eq:rhoFixed}
\rho_m=(0,1,3,5,\ldots,2m-1),
\qquad m\ge0.
\end{equation}
For example, the elementary part of $\psi$ is as follows: if the smallest part of $\rho$ is positive, one inserts a zero part; if the smallest part is $0$ and the next smallest part is even, one deletes the zero part.  The remaining case, in which the smallest part is $0$ and the next smallest part is odd, is handled by the modular hook map of Chen and Liu. We define the composition map
\begin{equation}\label{eq:iotaTransport}
\iota=\Phi^{-1}\circ\psi\circ\Phi.
\end{equation}
Noting that $\Phi$ is a bijection, we immediately see that $\iota$ is a complete weight-preserving, sign-reversing involution on all triples except the preimages of \eqref{eq:rhoFixed}.  These preimages are exactly
\begin{equation}\label{eq:tripleFixed}
X_m=(0;\varnothing,(1,2,\ldots,m)).
\end{equation}
They have
\[
w(X_m)=1+2+\cdots+m=\frac{m(m+1)}2,
\qquad
\operatorname{sgn}(X_m)=(-1)^m.
\]
Consequently
\begin{equation}\label{eq:FthetaNew}
F(q)=\sum_{m\ge0}(-1)^m q^{m(m+1)/2}.
\end{equation}
Combining \eqref{eq:S1tosignednew} and \eqref{eq:FthetaNew} proves the result.  An independent analytic derivation of \eqref{eq:FthetaNew}, obtained by specializing the refinement $\mathcal V(z;q)$ at $z=-1$, is given in Remark~\ref{rem:specializations} below.

\section{Moments of the length-difference statistic}\label{s:moments}
In this section we study the refinement $\mathcal V(z;q)$ defined in \eqref{eq:Vdef} and its moments \eqref{eq:momentdef}.

\subsection{Notation}
We write
\[
 D_z:=z\frac{\partial}{\partial z},
 \qquad
 \delta:=q\frac{\partial}{\partial q},
\]
and recall that $\Theta(q)=\sum_{r\ge0}(-1)^rq^{r(r+1)/2}$ from \eqref{eq:ABcompanion}.  Let
\begin{equation}
\label{eq:BPPsi}
\begin{aligned}
 B(z;q)&:=(1+z)(-zq,-z^{-1}q;q)_\infty,\\ P(z;q)&:=\sum_{r\geq0}z^{3r}q^{r(3r+1)/2}
          \bigl(1-z^2q^{2r+1}\bigr),\\
 \Psi(z;q)&:=\sum_{r\geq0}(-1)^r z^{2r+1}q^{r(r+1)/2},
\end{aligned}
\end{equation}
so that $\Psi(1;q)=\Theta(q)$.  Set
\begin{equation}
\label{eq:basicoperators}
 Q(q):=(q;q)_\infty,
 \qquad
 \epsilon(q):=\frac{\delta Q(q)}{Q(q)},
\end{equation}
\[
 B_0(q):=B(1;q),
 \qquad
 P_0(q):=P(1;q),
\]
and define the $q$-differential operators
\begin{equation}
\label{eq:operators}
 \mathscr L_B:=2\delta+2\epsilon+\frac14,
 \qquad
 \mathscr L_P:=6\delta+\frac14,
 \qquad
 \mathscr L_\Psi:=8\delta+1.
\end{equation}

\subsection{A partial theta decomposition}
The starting point is Ramanujan's partial theta identity \cite{AndrewsBerndt} (see also \cite{Alladi,Warnaar})
\begin{align}
\label{eq:Ramanujan}
 \sum_{n\geq0}\frac{q^n}{(xq;q)_n(q/x;q)_n}
 &=(1-x)\sum_{r\geq0}(-1)^r x^{3r}q^{r(3r+1)/2}
   \bigl(1-x^2q^{2r+1}\bigr)\nonumber\\
 &\quad+\frac{1}{(xq,q/x;q)_\infty}
 \sum_{r\geq0}(-1)^r x^{2r+1}q^{r(r+1)/2}.
\end{align}

\begin{lemma}
\label{lem:decomposition}
We have
\begin{equation}
\label{eq:Vdecomp}
 \mathcal{V}(z;q)=B(z;q)P(z;q)-\Psi(z;q).
\end{equation}
\end{lemma}

\begin{proof}
From \eqref{eq:Vdef},
\[
 \mathcal{V}(z;q)=(-zq,-z^{-1}q;q)_\infty
 \sum_{a\geq0}\frac{q^a}{(-zq;q)_a(-z^{-1}q;q)_a}.
\]
Set \(x=-z\) in \eqref{eq:Ramanujan} and multiply by $(-zq,-z^{-1}q;q)_\infty$.  Since
\[
 (-1)^r(-z)^{3r}=z^{3r},
 \quad
 (-1)^r(-z)^{2r+1}=-(-1)^rz^{2r+1},
\]
the first term on the right of \eqref{eq:Ramanujan} becomes $B(z;q)P(z;q)$ and the second becomes $-\Psi(z;q)$.
\end{proof}

\begin{remark}
\label{rem:specializations}
The specializations $z=\pm1$ of \eqref{eq:Vdecomp} connect with Sections~\ref{s1} and \ref{s:combnew}.
\begin{enumerate}
\item At $z=1$ we get $S_1=B_0P_0-\Theta$, that is,
\begin{equation}
\label{eq:P0}
 P_0(q)=\frac{S_1(q)+\Theta(q)}{B_0(q)}.
\end{equation}
Since $B_0=2(-q;q)_\infty^2=2(q^2;q^2)_\infty^2/(q;q)_\infty^2$, combining this with \eqref{eq:ABcompanion} gives
\[
 2q(q;q)_\infty T_o(q)=\frac{(q^2;q^2)_\infty^2}{(q;q)_\infty}P_0(q)-(q;q)_\infty\Theta(q),
\]
which is the partial theta identity for $T_o(q)$ obtained by Andrews and El Bachraoui \cite[(4.8)]{AB}.
\item At $z=-1$, the factor $1+z$ gives $B(-1;q)=0$, while $\Psi(-1;q)=-\Theta(q)$.  Hence
\[
 F(q)=\mathcal V(-1;q)=\Theta(q),
\]
which is \eqref{eq:FthetaNew}.  Together with \eqref{eq:S1tosignednew}, this gives an analytic counterpart of the involutive proof of Theorem~\ref{thm1AB}.
\end{enumerate}
\end{remark}

\subsection{Heat equations}
Each of the three pieces in \eqref{eq:Vdecomp} is annihilated by a heat-type operator.

\begin{lemma}
\label{lem:heat}
We have
\begin{align}
 (D_z^2+D_z-6\delta)P&=0,\label{eq:heatP}\\
 (D_z^2-8\delta-1)\Psi&=0,\label{eq:heatPsi}\\
 (D_z^2-D_z-2\delta-2\epsilon)B&=0.\label{eq:heatB}
\end{align}
\end{lemma}

\begin{proof}
In the first part of \(P\), the $z$- and \(q\)-exponents are, respectively,
\[
 m=3r,
 \qquad
 n=\frac{r(3r+1)}2,
\]
and they satisfy \(m^2+m=6n\).  In the second part they are
\[
 m=3r+2,
 \qquad
 n=\frac{r(3r+1)}2+2r+1,
\]
and again \(m^2+m=6n\).  Since $D_z$ and $\delta$ act on $z^mq^n$ by multiplication by $m$ and $n$, this proves \eqref{eq:heatP}.  For \(\Psi\), the exponents are \(m=2r+1\) and \(n=r(r+1)/2\), so \(m^2=8n+1\), which gives \eqref{eq:heatPsi}.  Finally, the Jacobi triple product \cite{AndrewsTOP} gives
\begin{equation}
\label{eq:JTP}
Q(q)B(z;q)=\sum_{m\in\mathbb{Z}}z^m q^{m(m-1)/2}.
\end{equation}
Each monomial on the right satisfies \(m^2-m=2\cdot m(m-1)/2\).  Thus
\[
 (D_z^2-D_z-2\delta)(QB)=0.
\]
Since \(\delta Q=\epsilon Q\) and $Q$ does not depend on $z$, division by \(Q\) gives \eqref{eq:heatB}.
\end{proof}

\subsection{The moment formula}
\begin{theorem}
\label{thm:moment}
For every \(k\geq0\),
\begin{equation}
\label{eq:allmoment}
 \mathcal{M}_{2k}(q)
 =\sum_{j=0}^{k}\binom{2k}{2j}
 \bigl(\mathscr L_B^{\,j}B_0\bigr)
 \bigl(\mathscr L_P^{\,k-j}P_0\bigr)
 -\mathscr L_\Psi^{\,k}\Theta.
\end{equation}
In particular,
\begin{equation}
\label{eq:M2}
 \mathcal{M}_2
 =6\delta S_1+\mathcal{E} S_1-2\delta\Theta+(\mathcal{E}-1)\Theta,
 \qquad
 \mathcal{E}=\frac{5E_2(\tau)-8E_2(2\tau)+9}{12},
\end{equation}
where \(q=e^{2\pi i\tau}\) and
\[
 E_2(\tau)=1-24\sum_{n\geq1}\sigma_1(n)q^n.
\]
\end{theorem}

\begin{proof}
Put \(z=e^t\) and define
\begin{equation}
\label{eq:tilde}
 \widetilde B(t;q):=e^{-t/2}B(e^t;q),
 \quad
 \widetilde P(t;q):=e^{t/2}P(e^t;q).
\end{equation}
The exponential factors cancel, so
\[
 B(e^t;q)P(e^t;q)=\widetilde B(t;q)\widetilde P(t;q).
\]
Also, \(B(z^{-1};q)=z^{-1}B(z;q)\), and therefore \(\widetilde B(-t;q)=\widetilde B(t;q)\).  In particular, every odd derivative of \(\widetilde B\) at \(t=0\) is zero. Because \(D_z\) becomes \(\partial_t\), Lemma~\ref{lem:heat} becomes
\begin{align}
 \partial_t^2\widetilde B&=\mathscr L_B\widetilde B,\label{eq:tildeheatB}\\
 \partial_t^2\widetilde P&=\mathscr L_P\widetilde P,\label{eq:tildeheatP}\\
 \partial_t^2\Psi(e^t;q)&=\mathscr L_\Psi\Psi(e^t;q).\label{eq:tildeheatPsi}
\end{align}
For instance, writing $B(e^t;q)=e^{t/2}\widetilde B$, we have $(\partial_t^2-\partial_t)B(e^t;q)=e^{t/2}\bigl(\partial_t^2-\tfrac14\bigr)\widetilde B$, so \eqref{eq:heatB} is equivalent to \eqref{eq:tildeheatB}; similarly $(\partial_t^2+\partial_t)P(e^t;q)=e^{-t/2}\bigl(\partial_t^2-\tfrac14\bigr)\widetilde P$ shows that \eqref{eq:heatP} is equivalent to \eqref{eq:tildeheatP}.  The operators on the right of \eqref{eq:tildeheatB}--\eqref{eq:tildeheatPsi} involve only \(q\), so they commute with \(\partial_t\).  By iteration, and since $\widetilde B(0;q)=B_0$, $\widetilde P(0;q)=P_0$, and $\Psi(1;q)=\Theta$,
\begin{align}
 \left.\partial_t^{2j}\widetilde B(t;q)\right|_{t=0}
 &=\mathscr L_B^{\,j}B_0,\label{eq:iteratedB}\\
 \left.\partial_t^{2j}\widetilde P(t;q)\right|_{t=0}
 &=\mathscr L_P^{\,j}P_0,\label{eq:iteratedP}\\
 \left.\partial_t^{2j}\Psi(e^t;q)\right|_{t=0}
 &=\mathscr L_\Psi^{\,j}\Theta.\label{eq:iteratedPsi}
\end{align}
On the other hand, by \eqref{eq:Vdef} and \eqref{eq:momentdef},
\[
 \mathcal{M}_{2k}(q)
=\left.\partial_t^{2k}\mathcal{V}(e^t;q)\right|_{t=0}.
\]
Using \eqref{eq:Vdecomp}, applying Leibniz's rule to
\(\widetilde B\widetilde P\), and using the vanishing of the odd derivatives of \(\widetilde B\), we get
\begin{align*}
 \mathcal{M}_{2k}(q)
 &=\sum_{j=0}^{k}\binom{2k}{2j}
 \left.\partial_t^{2j}\widetilde B\right|_{t=0}
 \left.\partial_t^{2k-2j}\widetilde P\right|_{t=0}
 -\left.\partial_t^{2k}\Psi(e^t;q)\right|_{t=0}.
\end{align*}
Substitution of \eqref{eq:iteratedB}--\eqref{eq:iteratedPsi} proves \eqref{eq:allmoment}.

It remains to simplify the case \(k=1\).  Let
\[
 b(q):=\frac{\delta B_0}{B_0}.
\]
Since \(B_0=2(-q;q)_\infty^2=2(q^2;q^2)_\infty^2/(q;q)_\infty^2\), standard logarithmic differentiation gives
\begin{equation}
\label{eq:bandeps}
 b=\frac{2E_2(2\tau)-E_2(\tau)-1}{12},
 \qquad
 \epsilon=\frac{E_2(\tau)-1}{24}.
\end{equation}
Formula \eqref{eq:allmoment} with \(k=1\) is
\[
 \mathcal{M}_2
 =B_0\mathscr L_P P_0
 +(\mathscr L_B B_0)P_0
 -\mathscr L_\Psi\Theta.
\]
Using \eqref{eq:P0}, we find
\begin{align*}
 B_0\mathscr L_P P_0
 &=6\delta(S_1+\Theta)+\left(\frac14-6b\right)(S_1+\Theta),\\
 (\mathscr L_B B_0)P_0
 &=\left(2b+2\epsilon+\frac14\right)(S_1+\Theta),\\
 \mathscr L_\Psi\Theta&=8\delta\Theta+\Theta.
\end{align*}
Therefore
\[
 \mathcal{M}_2
 =6\delta S_1-2\delta\Theta
 +\left(-4b+2\epsilon+\frac12\right)(S_1+\Theta)-\Theta.
\]
By \eqref{eq:bandeps},
\[
 -4b+2\epsilon+\frac12
 =\frac{5E_2(\tau)-8E_2(2\tau)+9}{12}=\mathcal{E}.
\]
This is exactly \eqref{eq:M2}.
\end{proof}

\begin{example}
The first few values of $M_2(n)$ are
\[
 \mathcal M_2(q)=2q+2q^2+12q^3+14q^4+26q^5+56q^6+78q^7+120q^8+182q^9+\cdots,
\]
in agreement with \eqref{eq:M2}.
\end{example}

\subsection{Proof of Theorem~\ref{thm:cong}}
Multiply \eqref{eq:Vdecomp} by \(Q(q)=(q;q)_\infty\).  From \eqref{eq:JTP} and Euler's pentagonal number theorem
\[
 Q(q)=\sum_{s\in\mathbb{Z}}(-1)^s q^{s(3s-1)/2},
\]
we obtain, after applying \(D_z^{2k}\) and setting \(z=1\),
\begin{align}
\label{eq:representation}
 d_{2k}(n)&=
 \sum_{\substack{r\geq0,\ m\in\mathbb{Z}\\
 \frac{m(m-1)}2+\frac{r(3r+1)}2=n}}
 (m+3r)^{2k}-\sum_{\substack{r\geq0,\ m\in\mathbb{Z}\\
 \frac{m(m-1)}2+\frac{r(3r+1)}2+2r+1=n}}
 (m+3r+2)^{2k}\nonumber\\
 &\hspace{1cm}-\sum_{\substack{r\geq0,\ s\in\mathbb{Z}\\
 \frac{r(r+1)}2+\frac{s(3s-1)}2=n}}
 (-1)^{r+s}(2r+1)^{2k}.
\end{align}
This holds for every $k\ge0$ (with $0^0=1$).  Every term in these three sums comes from a representation of the same integer by \(x^2+3y^2\).  In the first sum put
\[
 x=6r+1,
 \qquad
 y=2m-1.
\]
Then
\[
 x^2+3y^2=4(6n+1),
 \qquad
 m+3r=\frac{x+y}{2}.
\]
In the second sum put
\[
 x=6r+5,
 \qquad
 y=2m-1.
\]
Again,
\[
 x^2+3y^2=4(6n+1),
 \qquad
 m+3r+2=\frac{x+y}{2}.
\]
In the third sum put
\[
 x=6s-1,
 \qquad
 y=2r+1.
\]
Then
\[
 x^2+3y^2=4(6n+1),
 \qquad
 2r+1=y.
\]
Thus every summand in \eqref{eq:representation} is a signed \(2k\)-th power of either \((x+y)/2\) or \(y\), where
\begin{equation}
\label{eq:norm}
 x^2+3y^2=4(6n+1).
\end{equation}
Now let \(p\equiv5\pmod6\).  Then \(-3\) is a quadratic nonresidue modulo \(p\).  Hence
\begin{equation}
\label{eq:inert}
 x^2+3y^2\equiv0\pmod p
 \quad\Longrightarrow\quad
 x\equiv y\equiv0\pmod p.
\end{equation}
Indeed, if \(y\not\equiv0\pmod p\), then \((xy^{-1})^2\equiv-3\pmod p\), a contradiction; if \(y\equiv0\), then also \(x\equiv0\).

Suppose first that \(v_p(6n+1)\) is odd.  If \eqref{eq:norm} had a solution, then \eqref{eq:inert} would imply \(p\mid x,y\).  Dividing by \(p^2\) and repeating shows that the \(p\)-adic valuation of the represented norm must be even.  This is impossible.  Therefore \eqref{eq:norm} has no integral solutions, all three sums in \eqref{eq:representation} are empty, and \(d_{2k}(n)=0\).

Suppose next that \(v_p(6n+1)=2a\).  Repeated use of \eqref{eq:inert} gives
\[
 p^a\mid x,
 \qquad
 p^a\mid y
\]
for every solution of \eqref{eq:norm}.  Since \(p\) is odd, both \((x+y)/2\) and \(y\) are divisible by \(p^a\).  Every individual weight in \eqref{eq:representation} is therefore divisible by \(p^{2ka}\), and the same is true of their signed sum.  Hence
\[
 p^{2ka}\mid d_{2k}(n).
\]
Finally, put
\[
 N=pn+\frac{5p-1}{6},
\]
which is an integer because $5p\equiv1\pmod 6$.  Then
\[
 6N+1=p(6n+5).
\]
If its \(p\)-adic valuation is odd, then \(d_{2k}(N)=0\).  If it is even, it is at least \(2\), so the preceding divisibility gives \(p^{2k}\mid d_{2k}(N)\).  This proves \eqref{eq:progression}.\qed

\begin{example}
The first instances of \eqref{eq:progression} are
\[
 d_2(5n+4)\equiv0\pmod{25},
 \qquad
 d_4(5n+4)\equiv0\pmod{625},
\]
and
\[
 d_2(11n+9)\equiv0\pmod{121}.
\]
\end{example}

\begin{remark}
The important point is that these congruences are termwise.  The prime \(p\equiv5\pmod6\) divides both coordinates of every contributing representation, and the \(2k\)-th moment turns this into divisibility by \(p^{2k}\).
\end{remark}

The case $k=0$ of the vanishing argument gives information about the companion $C(q)$ itself.

\begin{corollary}
\label{cor:cvanish}
If $p\equiv5\pmod6$ is prime and $v_p(6n+7)$ is odd, then $c(n)=0$.
\end{corollary}

\begin{proof}
Since $\mathcal M_0=S_1$, multiplying \eqref{eq:ABcompanion} by $(q;q)_\infty$ gives
\[
 \sum_{n\ge0}d_0(n)q^n=(q;q)_\infty\Theta(q)+4qC(q).
\]
The coefficient of $q^{n+1}$ in $(q;q)_\infty\Theta(q)$ is exactly the third sum in \eqref{eq:representation} with $k=0$ and $n$ replaced by $n+1$.  If $v_p(6(n+1)+1)=v_p(6n+7)$ is odd, then \eqref{eq:norm} has no solutions, so $d_0(n+1)$ and this coefficient both vanish.  Hence $4c(n)=0$.
\end{proof}

For $p=5$ and $n=25m+j$ with $j\in\{8,13,18,23\}$, one has $6n+7=5(30m+11),\,5(30m+17),\,5(30m+23),\,5(30m+29)$, respectively, so $v_5(6n+7)=1$.  Thus Corollary~\ref{cor:cvanish} contains the vanishing statements of Theorem~\ref{thm2AB}.  This is consistent with \eqref{RD}: the integer $D=6n+7$ there is precisely $6(n+1)+1$.

\subsection{Special values}
Besides $z=\pm1$, the refinement $\mathcal V(z;q)$ has simple product-sum forms at roots of unity of small order.  For $z = i$, the factorization $(1+ix)(1-ix)=1+x^2$ gives
\[
\mathcal{V}(i;q)
=
\sum_{a=0}^{\infty}
q^a(-q^{2a+2};q^2)_{\infty}.
\]
For $z = \omega = e^{2\pi i/3}$, the factorization $(1+\omega x)(1+\omega^{-1}x)=(1+x^3)/(1+x)$ gives
\[
\mathcal{V}(\omega;q)
=
\frac{(-q^3;q^3)_\infty}{(-q;q)_\infty}
\sum_{a=0}^{\infty}
q^a
\frac{(-q;q)_a}{(-q^3;q^3)_a}.
\]
Since $N(m,n)=N(-m,n)$, both values are real, and the coefficients of $q^n$ in $\mathcal V(i;q)$ and $\mathcal V(\omega;q)$ are, respectively,
\[
 \sum_{m\equiv0\Mod{4}}N(m,n)-\sum_{m\equiv2\Mod{4}}N(m,n)
\]
and
\[
 \sum_{m\equiv0\Mod{3}}N(m,n)-\frac12\sum_{m\not\equiv0\Mod{3}}N(m,n).
\]

\section{Conclusion}\label{s:conclusion}
We conclude with a few open questions below.
\vspace{0.2cm}

\begin{enumerate}
    \item Is there a rank or crank statistic on the signed objects
associated with \(S_1(q)\) that explains the modulo \(4\) congruence directly?  The length difference $\ell(\lambda)-\ell(\mu)$ underlying $\mathcal V(z;q)$ is a natural candidate to start from.
\vspace{0.1cm}

\item Prove the Andrews--El Bachraoui modulo $8$ congruence \cite[Conjecture 2]{AB} for $s_1(n)$ combinatorially.
\vspace{0.1cm}

\item Do the special values $\mathcal V(i;q)$ and $\mathcal V(\omega;q)$ admit partial theta decompositions analogous to \eqref{eq:Vdecomp}, and do they lead to congruences or sign patterns for the corresponding residue-class counts of $N(m,n)$?
\end{enumerate}

\end{document}